\documentclass{article}

\usepackage{amsthm}
\usepackage[utf8]{inputenc}
\usepackage[T1]{fontenc}
\usepackage[centertags]{amsmath}
\usepackage{amsfonts}
\usepackage{graphicx}
\usepackage{amsmath}
\usepackage{amssymb}
\usepackage{latexsym}
\usepackage[mathscr]{euscript}

\usepackage{tikz,color,makeidx} 
\usetikzlibrary{matrix,arrows,chains,positioning,scopes,decorations.pathmorphing,backgrounds,fit}

\numberwithin{equation}{section}
\newtheorem{teo}{Theorem}[section]
\newtheorem{cor}{Corollary}[section]
\newtheorem{pro}{Proposition}[section]
\newtheorem{lem}{Lemma}[section]
\newtheorem{conj}{Conjecture}[section]
\newtheorem{es}{\textbf{Example}}[section]
\newtheorem{defi}{\mbox{\textbf{Definition}}}[section]
\newtheorem{rem}{\mbox{\textbf{Remark}}}[section]

\newcommand{\bdfn}{\begin{defi} \begin{rm}}
\newcommand{\edfn}{\end{rm} \end{defi}}
\newcommand{\bthm}{\begin{teo}}
\newcommand{\ethm}{\end{teo}}
\newcommand{\bprop}{\begin{pro}}
\newcommand{\eprop}{\end{pro}}
\newcommand{\bcor}{\begin{cor}}
\newcommand{\ecor}{\end{cor}}
\newcommand{\blem}{\begin{lem}}
\newcommand{\elem}{\end{lem}}
\newcommand{\bfact}{\begin{rem} \begin{rm}}
\newcommand{\efact}{\end{rm} \end{rem}}
\newcommand{\bex}{\begin{es} \begin{rm}}
\newcommand{\eex}{ \end{rm} \end{es}}

\newcommand{\vsp}[1]{\vspace*{#1cm}}

\newcommand{\quot}[2]{{\raisebox{.2em}{$#1$}\left/\raisebox{-.2em}{$#2$}\right.}}

\tikzset{node distance=2cm, auto}

\title{Towards understanding the Pierce-Birkhoff conjecture via MV-algebras}

\author{Serafina Lapenta\\ 
 {\small Department of Mathematics, Computer Science and Economics,}\\ 
{\small University of Basilicata,
 Viale dell'Ateneo Lucano, 10 Potenza, Italy}\\
 {\small serafina.lapenta@unibas.it}
\and 
Ioana Leu\c stean \\
{\small Department of Computer Science,} \\
{\small Faculty of Mathematics and Computer Science, University of Bucharest,}\\
{\small Academiei nr.14, sector 1, C.P. 010014,  Bucharest, Romania}\\   
{\small ioana@fmi.unibuc.ro}
}

\date{}

\begin{document}
\maketitle 

\begin{abstract}
Our main issue was to understand the connection  between \L ukasiewicz logic with product and the Pierce-Birkhoff conjecture, and to express it in a  mathematical way.  To do this we define the class of \textit{f}MV-algebras, which are MV-algebras endowed with both an internal binary product and a scalar product with scalars from $[0,1]$.
The proper quasi-variety  generated by $[0,1]$,  with both products interpreted as the real product,  provides the desired framework:
the normal form theorem of its corresponding logical system can be seen as a local version of the Pierce-Birkhoff conjecture.

\end{abstract}

\section*{Introduction}

\noindent \L ukasiewicz  $\infty$-valued logic with the  primitive connectives  $\rightarrow$ and $\neg$ has as set of truth values the real interval $[0,1]$ with $x\rightarrow y=\max(0,1-x+y)$ and $\neg x=1-x$ for any $x,y\in [0,1]$. The completeness theorem states that a formula is provable if and only if it holds in the standard model $[0,1]$.
The corresponding algebraic structures, called {\em MV-algebras} \cite{Cha1}, were defined as  structures $(A,\oplus,^*,0)$ of type $(2,1,0)$ satisfying some appropriate axioms. The variety of MV-algebras is generated by $([0,1], \oplus, ^*,0)$ where 
$x\oplus y=\min(1,x+y)$ and  $x^*=1-x$ for any $x, y\in [0,1]$. 

One major result in the theory of MV-algebras is the categorical equivalence with the Abelian  lattice-ordered groups with strong unit \cite{Mun1}. As a consequence, for any MV-algebra $A$ there exists an Abelian  lattice-ordered group with strong unit $(G,u)$ such that $A\simeq [0,u]_G$, where $[0,u_G]=([0,u], \oplus, ^*,0)$ with $x\oplus y= u\wedge (x+y)$ and  $x^*=u-x$ for any $x, y\in [0,u]$.

 The standard MV-algebra $[0,1]$ is closed with respect to the real product, so  the study of the structure obtained  by  endowing the MV-algebras with a  product operation was an active research direction. Our  model in this case is $([0,1], \oplus, \cdot, ^*, 0)$, 
where   $([0,1], \oplus, ^*, 0)$ is the standard MV-algebra and $\cdot$ is the real product.
We shall only mention the references strictly related with our development, but the literature on the subject is far richer.  For such structures further categorical equivalences with 
classes of lattice-ordered structures were proved. 

If we consider the real product as a binary operation on $[0,1]$, the corresponding structures are MV-algebras $A$ endowed 
with  an operation $\cdot:A\times A\to A$. These were introduced and studied in \cite{DiND}  under the name of {\em PMV-algebras} and they are categorically equivalent to a class of lattice-ordered rings with strong unit.  Particular important subclasses were further studied in \cite{MonPMV, Mon+}. As proved in \cite{MonPMV}, the variety of (commutative and unital) PMV-algebras is larger than the variety generated by $([0,1], \oplus,\cdot, ^*, 0)$. The quasi-variety 
$ISP([0,1])$ was characterized in \cite{Mon+}. 

One can also consider the real product on $[0,1]$ as a scalar multiplication. Our standard  model in this case 
is $([0,1],\oplus, ^*, \{\alpha\mid \alpha\in [0,1]\}, 0)$ where $x\mapsto \alpha x$ is a unary operation for any $\alpha\in [0,1]$. 
These structures are investigated in \cite{LeuRMV} under the name of {\em Riesz MV-algebras} and they are categorically equivalent to Riesz spaces (vector lattices)  with strong unit. The variety of Riesz MV-algebras is generated by the standard model $[0,1]$. 

Both for PMV-algebras and Riesz MV-algebras logical systems were developed \cite{HC, LeuRMV} and such systems are conservative extensions of \L ukasiewicz logic.  One of the main theorems of \L ukasiewicz logic states that the term functions corresponding to the formulas of \L ukasiewicz logic with $n$ variables ($n\geq 1$) are exactly the continuous $[0,1]$-valued piecewise linear functions 
with integer coefficients defined on $[0,1]^n$ \cite{McN}. This  can  be seen as a {\em normal form theorem} for \L ukasiewicz logic. A similar result was proved in \cite{LeuRMV} for the logical system that has Riesz MV-algebras as models; in this case the piecewise linear functions have real coefficients.  In \cite[Introduction]{MoPa} it is stated that a similar result for PMV-algebras   is related to the Pierce-Birkhoff conjecture \cite{BP} (see Section \ref{pbc} for explicit formulations). 

Our main issue was to understand the connection  between \L ukasiewicz logic with product (in a broad sense) and the Pierce-Birkhoff conjecture, and to express it in a mathematical way.  To do this,  we study the class of structures obtained by endowing MV-algebras with both the internal binary product and the scalar product  (as a family of unary operations).  We called {\em $\textit{f}$MV-algebras} the structures obtained in this way and the model 
we focus on is ${\mathbf I}=([0,1],\oplus,\cdot, ^*, \{\alpha\mid \alpha\in [0,1]\}, 0)$.

We  can briefly express our final  conclusions as follows: \\
(1) the variety of \textit{f}MV-algebras is larger then HSP$({\mathbf I})$;\\
(2) we characterize ISP$({\mathbf I})$; we called  {\em FR$^+$-algebras} the members of ISP$({\mathbf I})$;\\  
(3) the logical system  $\mathcal{FMVL}^{+}$ that has 
FR$^+$-algebras as models is a conservative extension of \L ukasiewicz logic and it is  complete w.r.t. $\mathbf I$; \\
(4) the normal form theorem  for $\mathcal{FMVL}^{+}$  (Conjecture \ref{end}) is, in our approach,   the link between the  Pierce-Birkhoff conjecture and \L ukasiewicz logic with product. We proved it for $n\leq 2$, due to the fact that the Pierce-Birkhoff conjecture is proved in this case \cite{Mahe}. Note that Conjecture \ref{end} does not immediately imply, nor it is implied by the Pierce-Birkhoff conjecture, 
additional results are needed and, in our opinion, these results belong  more to the area of algebraic geometry than to the area of logic.  

Note that for FR$^+$-algebras we prove  the subdirect representation w.r.t. totally-ordered structures and  a representation
as algebras of $^*[0,1]$-valued functions, where $^*[0,1]$ is an ultrapower of $[0,1]$.

Section \ref{sec0} contains the preliminary notions.  In Section \ref{sec1} we define the \textit{f}MV-algebras, we prove a subdirect representation theorem and  the categorical equivalence between  \textit{f}MV-algebras and a subclass of \textit{f}-algebras. In Section \ref{sec2} we investigate special classes of \textit{f}MV-algebras and in the last subsection we introduce the class of FR$^+$-algebras. Logical systems are developed in Section \ref{sec3}.  In Section \ref{pbc} we relate  previous results with the Pierce-Birkhoff conjecture.

\section{Preliminaries}\label{sec0}

\noindent We  recall the algebraic structures involved in our development. 

\subsection{MV-algebras and $\ell u$-groups}
\noindent An {\em MV-algebra} is an algebraic structure $(A,\oplus, ^*, 0)$  such that $(A, \oplus, 0)$ is an abelian monoid, $(x^*)^*=x$ and $(x^{*}\oplus y)^{*}\oplus y=(y^{*}\oplus x)^{*}\oplus x$ for any $x,y\in A$.  If $A$ is an MV-algebra, one can define $x\odot y= (x^*\oplus  y^*)^*,$
 $x\ominus y=x\odot y^*,$ $x\rightarrow y=x^* \oplus y$  and $1=0^*$, for any $x,y\in A$.  The  order defined by setting $x\le y$ if and only if $x\odot  y^*=0$ is a  lattice order such that 
  $  x\vee y= (x\odot  y^*)\oplus y$ and $
  x\wedge y= (  x^* \vee y^*)^*=x\odot ( x^*\oplus y)$ for any $x,y\in A$.
We refer to \cite{CDM} for all the unexplained notions concerning MV-algebras and to \cite{DM} for advanced topics. 

The variety of MV-algebras is generated by  $([0,1], \oplus, ^*,0)$ where 
$x\oplus y=\min(1,x+y)$ and  $x^*=1-x$ for any $x, y\in [0,1]$. Any MV-algebra is a subdirect product of totally-ordered MV-algebras. Moreover, any MV-algebra can be embedded into a subalgebra of an ultrapower of $[0,1]$.

An ideal $I$ in an  MV-algebra $A$ is a  nonempty subset such that for any $x,y\in I$, $x\oplus y\in I$ and if $x\le y$ with $y\in I$ then $x\in I$. The intersection of all maximals ideals in an MV-algebra is called \textit{radical}, and it is denoted by $Rad(A)$. An  MV-algebra is called \textit{semisimple} if $Rad(A)=\{ 0 \}$. Recall that  any semisimple MV-algebra is isomorphic to a subalgebra of the MV-algebra $C(X)$ of $[0,1]$-valued continuous functions defined on a suitable topological space $X$.

The relation between MV-algebras and lattice-ordered groups plays a crucial role in our approach.
A {\em  lattice-ordered group}  is a group is a structure \mbox{$(G,+,0,\leq)$} such that $(G,+,0)$ is a group, $(G,\leq)$ is a lattice   and $x\le y$ implies $x+z\le y+z$ for any $x,y,z\in G$ \cite{BigKim,Birk}. A strong unit $u$  is an  element $u>0$ such that for any $x\in G$ there exists an $n\in \mathbb{N}$ with  $x\le nu$.  An {\em $\ell  u$-group} is a pair $(G,u)$ where $G$ is an Abelian lattice-ordered group and $u$ is a strong unit in $G$.
For any element $u\in G$, $u>0$ we define the \textit{unit interval} as
\vsp{-0.2}
\begin{center}
$ [0,u]_G=([0,u]=\{ x\in G \mid 0\le x\le u\},\oplus, ^*, 0)$, where
\end{center}
\vsp{-0.2}
$ x\oplus y=u\wedge (x+y)$ and $x^*=u-x$ for any $x,y\in [0,u]$.
The structure $[0,u]_G$ is an MV-algebra.

An $\ell$-group $G$ is said to be {\em Archimedean} if $nx\le y$ for any $n\in \mathbb{N}$, with $x,y\in G$, implies $x\le 0$.

Let $\mathbf{MV}$ denote the category of MV-algebras and 
 $\mathbf{auG}$ the category whose objects are $\ell u$-groups and whose morphisms are  homomorphisms of $\ell u$-groups that preserve the strong unit.
\bthm \cite{Mun1}
Let  $\Gamma$ be the functor defined by 
\vsp{-0.2}
\begin{center}
$\Gamma: \mathbf{auG}\rightarrow \mathbf{MV}$ \quad $\Gamma (G,u)=[0,u]_G $ and $\Gamma(h)= h|_{[0,u]_G},$
\end{center}
\vsp{-0.2}
then $\Gamma$ establishes a categorical equivalence. Moreover,
an MV-algebra is semisimple if and only if the corresponding $\ell u$-group is Archimedean.
\ethm

\subsection{PMV-algebras and $\ell u$-rings} \label{subsPMV}

\noindent A  {\em  PMV-algebra} \cite{DiND} is a  structure $(P, \oplus, \ ^*, \cdot , 0)$ such that $(P, \oplus, \ ^*, 0)$ is an MV-algebra and the binary operation $\cdot $ satisfies the following, for any $x,y,z \in P$:

(PMV1) $z\cdot (x\odot (x\wedge y)^*)=(z\cdot x)\odot (z\cdot (x\wedge y))^*$

(PMV2) $(x\odot (x\wedge y)^*) \cdot z=(x\cdot z)\odot ((x\wedge y)\cdot z)^*$.

(PMV3) $x \cdot (y\cdot z)= (x \cdot y)\cdot z$.\\
A PMV-algebra that has unit for product is called \textit{unital}, and in this case  the unit is $1=0^*$  \cite[Proposition 3.1]{DiND}. 

A {\em $\cdot$-ideal} $I$ in a PMV-algebra $P$ is an MV-ideal that satisfies the condition $x\cdot y \in I$ and $y\cdot x \in I$ for any $x\in I$ and $y\in P$.

A {\em PMV\textit{f}-algebra} is an PMV-algebra that satisfies the following condition:

(f) $x\wedge y=0$ implies $(x\cdot z)\wedge y=(z\cdot x)\wedge y=0$ for any $x,y,z\in P$.\\
By \cite[Theorem 5.4]{DiND} this condition has an equational form. Trivially, unital PMV-algebras are PMV\textit{f}-algebras and 
the standard  MV-algebra [0,1] endowed with the usual real product is a  PMV\textit{f}-algebra. 
Any PMV\textit{f}-algebra is a subdirect product of totally-ordered PMV\textit{f}-algebra. 

The {\em Product MV-algebras} defined in \cite[Definition 2.8]{MonPMV} are the  unital and commutative PMV-algebras from \cite{DiND} (\cite{MonPMV}, \cite[Theorem 2.8]{Mon+}). A  $PMV^+$ -algebra $P$ is a unital and commutative PMV-algebra that satisfies the quasi-identity $(x^2=0 \Rightarrow x=0)$ for any $x\in P$ \cite{Mon+}. Any $PMV^+$-algebra is a subdirect product of totally-ordered $PMV^+$-algebras \cite{HC}.

 \bthm \cite{Mon+}
The class of all $PMV^+$-algebras  is the quasi-variety generated by $([0,1],\oplus, \cdot, ^*, 0)$. Moreover, any totally-ordered $PMV^+$-algebra embeds into a subalgebra of an ultrapower of $[0,1]$.
\ethm

An {\em $\ell u$-ring} $R$ is an $\ell u$-group endowed with an operation $\cdot :R \times R \rightarrow R$ such that $(R,\cdot)$ is a ring and   $x\cdot y\ge 0$ whenever $x, y\ge 0$ \cite{BigKim,BP}.  An {\em $f$-ring} is a lattice-ordered ring $R$ such that  $x\wedge y=0$ implies $(x\cdot z)\wedge y=(z\cdot x)\wedge y=0$, for any $z\geq 0$ and $x,y\in R$.

 Denoted by $\mathbf{PMV}$ the category of PMV-algebras and product preserving homomorphisms of MV-algebras and by $\mathbf{uR}$ the category of $\ell u$-rings and homomorphisms of $\ell u$-rings, we extend the functor $\Gamma$ and we get a functor $\Gamma_{(\cdot)}: \mathbf{uR}\rightarrow \mathbf{PMV}$, with $\Gamma_{(\cdot)}(R,u)=[0,u]_R $ and $\Gamma_{(\cdot)}(h)=h|_{[0,u]_R},$ where $(R,u)$ is an object and $h$ is a morphism in $\mathbf{uR}$. The functor $\Gamma_{(\cdot)}$ establishes a categorical equivalence \cite{DiND}. Moreover,  $f$-rings correspond to PMV\textit{f}-algebras.

\subsection{Riesz MV-algebras and  Riesz spaces}
\noindent A {\em Riesz MV-algebra} \cite{LeuRMV} is a structure $(R, \oplus, ^*, \{\alpha\mid \alpha\in [0,1]\},0)$ such that 
$(R, \oplus,  ^*, 0)$ is an MV-algebra and  $\{\alpha\mid \alpha\in [0,1]\}$ is a family of unary operations  satisfying the following identities for any $\alpha,\beta\in [0,1]$ and any $x,y\in R$:

(RMV1) $\alpha (x\odot y^*)= (\alpha x)\odot (\alpha y)^*$,

(RMV2) $\max(0,\alpha-\beta) x=(\alpha x)\odot (\beta x)^*$,

(RMV3) $\alpha  (\beta x)= (\alpha\beta) x$,

(RMV4) $1 x=x$.\\
By \cite[Corollary 2]{LeuRMV} any  homomorphism of MV-algebras between Riesz MV-algebras preserves the additional unary operations, so it is a homomorphism of Riesz MV-algebras.  Similarly, if $R$ is a Riesz MV-algebra then any ideal $I$ of the MV-algebra reduct is closed to the unary operations: $\alpha x\in I$  for any $x\in I$ and $\alpha\in [0,1]$.
The variety of Riesz MV-algebras is generated by $([0,1],\oplus, \cdot, ^*, \{\alpha\mid \alpha\in [0,1]\},  0)$ \cite{LeuRMV}.

A \textit{Riesz space} is a lattice ordered linear space over the real field, that is an $\ell$-group endowed with a scalar product with scalars over $\mathbb{R}$ such that $\alpha x\ge 0$ for any $\alpha$ and $x$ non-negative \cite{RSZan}. When we deal with $\ell u$-groups, we obtain a Riesz space with strong unit.

Let $\mathbf{uRS}$ be the category of Riesz spaces with strong unit and linear homomorphisms of $\ell u$-groups, and let $\mathbf{RMV}$ be the category of Riesz MV-algebras with homomorphisms of Riesz MV-algebras. We get a functor $\Gamma_{\mathbb{R}}: \mathbf{uRS}\rightarrow \mathbf{RMV}$, with $\Gamma_{\mathbb{R}}(R,u)=[0,u]_R $ and $\Gamma_{\mathbb{R}}(h)=h|_{[0,u]_R},$ where $(R,u)$ is an object and $h$ is a morphism in $\mathbf{uRS}$.
The functor $\Gamma_{\mathbb{R}}$ establishes a categorical equivalence \cite{LeuRMV}.

Putting all together, we get  the following commutative diagram

\begin{center}
\begin{tikzpicture}
  \node (A) {$\mathbf{uR}$};
  \node (B) [below of=A] {$\mathbf{PMV}$};
  \node (C) [right of=A] {$\mathbf{auG}$};
  \node (D) [below of=C] {$\mathbf{MV}$};
  \node (E) [right of=C] {$\mathbf{uRS}$};
  \node (F) [below of=E] {$\mathbf{RMV}$};
  \draw[->] (A) to node [swap] {$\Gamma_{(\cdot)}$} (B);
  \draw[->] (A) to node {$\mathcal{U}_{(\cdot \ell)}$} (C);
  \draw[->] (C) to node [swap] {$\Gamma$} (D);
  \draw[->] (B) to node [swap] {$\mathcal{U}_{(\cdot)}$} (D);
  \draw[->] (E) to node [swap] {$\mathcal{U}_{( \ell \mathbb{R})}$} (C);
  \draw[->] (E) to node {$\Gamma_{\mathbb{R}}$} (F);
  \draw[->] (F) to node {$\mathcal{U}_{\mathbb{R}}$} (D);
 \node [below=0.5cm, align=flush center] at (D){ Figure 1.};
\end{tikzpicture}
\end{center}
where $\mathcal{U}_{\mathbb{R}}$, $\mathcal{U}_{( \ell \mathbb{R})}$, $ \mathcal{U}_{(\cdot)}$, $ \mathcal{U}_{(\cdot \ell)}$ are  forgetful functors.\\

\subsection{ \textit{f}-algebras}

\noindent An  {\em \textit{f}-algebra} is a strucure $(L, +, \cdot, \{r|r\in{\mathbb R}\}, 0, \le )$   such that  $(L, +, \cdot, 0, \le)$ is a \textit{f}-ring, $(L, +, \{r|r\in{\mathbb R}\}, 0, \le)$ is a Riesz space and the following  condition is satisfied:

\vsp{-0.2}
\begin{center}
(fa) $r (x\cdot y)=(rx)\cdot y= x\cdot (ry)$ for any $x,y\in L$ and any $r\in \mathbb{R}$. 

\end{center}
\vsp{-0.2}
If there exists an element $e$ such that $xe=ex=x$ for any $x\in L$, $L$ will be called unital, and $e$ is the unit element for $L$. An \textit{f}-algebra with strong unit is an \textit{f}-algebra in which the underlying Riesz space is a Riesz space with strong unit. This notion was first introduced in \cite{BP}, for further details see also \cite{BigKim, HI, RSZan}.

\bdfn
In the sequel, an \textit{fu}{\em -algebra} is a pair $(V,u)$, where $V$ is an  \textit{f}-algebra and $u\in V$ is a strong unit  
such that $u\cdot u\leq u$.  If $(V_1,u_1)$ and $(V_2,u_2)$ are \textit{fu}-algebras then   $g:V_1\to V_2$ is a  homomorphism of \textit{fu}-algebras if it is a  homomorphism of \textit{f}-algebras  and $g(u_1)=u_2$.
\edfn

Note that, for an \textit{fu}-algebra $(V,u)$,   the interval $[0,u]$ is closed with respect to the product operations and with respect to the scalar multiplication by elements from $[0,1]$.

\section{\textit{f}MV-algebras}\label{sec1}

\subsection{Definitions and first results}
\noindent In the following we define the structure of \textit{fMV-algebra}, an appropriate notion of an ideal and we prove a subdirect representation theorem.

\bdfn\label{def:fMV}
An {\textit \textit{f}MV-algebra} is a structure $(A, \oplus,\cdot, ^*, \{\alpha\mid \alpha\in [0,1]\},0)$ which satisfies the following properties for any $x,y,z \in A$ and $\alpha\in [0,1]$:

(\textit{f}MV1) \quad  $(A,\oplus,\cdot,^*,0)$ is a PMV\textit{f}-algebra,

(\textit{f}MV2)\quad  $(A, \oplus, ^*, \{\alpha\mid \alpha\in [0,1]\},0)$  is a Riesz  MV-algebra,

(\textit{f}MV3)\quad $\alpha(x\cdot y)=(\alpha x)\cdot y= x\cdot (\alpha y)$.

\noindent We say that the \textit{f}MV-algebra  $A$ is \textit{unital}  (\textit{commutative}) if  its PMV-algebra reduct  $(A,\oplus,\cdot,^*,0)$ is unital (commutative).  
\edfn

\bfact
The class of   \textit{f}MV-algebras is a variety, denoted by $\mathbb{FMV}$.
\efact

One can see that an \textit{f}MV-algebra has an MV-algebra reduct, a Riesz MV-algebra reduct, as well as a PMV-algebra reduct and the forgetful functors are summarized in the following diagram:

\begin{center}
\begin{tikzpicture}
\node (C) {$\mathbf{fMV}$};
\node (D) [below of=C] {$\mathbf{MV}$};
\node (F) [right of=D] {$\mathbf{RMV}$};
\node (B) [left of=D] {$\mathbf{PMV}$};
\draw[->] (C) to node {$\mathcal{U}_{(\mathbb{R},\cdot )}$} (D);
\draw[->] (C) to node [swap] {$\mathcal{U}_{\mathbb{R}}$} (B);
\draw[->] (C) to node {$\mathcal{U}_{(\cdot)}$} (F);
\draw[->] (B) to node [swap] {$\mathcal{U}_{(\cdot)}$} (D);
\draw[->] (F) to node {$\mathcal{U}_{\mathbb{R}}$} (D);
\node [below=0.5cm, align=flush center] at (D){ Figure 2.};
\end{tikzpicture}
\end{center}

\blem \label{rem:01} If an  \textit{f}MV-algebra $A$ is unital then the unit is $1$.
\elem
\begin{proof}
If $A$ is a unital \textit{f}MV-algebra, then $\mathcal{U}_{\mathbb{R}}(A)$ is a unital PMV-algebra. The result follows from \cite[Proposition 3.1]{DiND}.
\end{proof}
 
The following examples play an important role in our development.
\bex \label{lem:01B}
\begin{enumerate}
\item 
The real interval $[0,1]$  with the standard structure of MV-algebra   is an \textit{f}MV-algebra  where both the scalar multiplication and the internal  product coincides with the  real product. We always consider $[0,1]$ endowed with this structure, but we note that it does not  generate the variety of all \textit{f}MV-algebras, but only a proper subclass (see Section \ref{formreal}).
\item If $(V,u)$ is an \textit{fu}-algebra then,   
 by \cite[Lemma 5]{LeuRMV} and \cite[Theorem 5.2]{DiND}, the interval $[0,u]$ is trivially a Riesz MV-algebra and a PMV\textit{f}-algebra. The condition (\textit{f}MV3) from Definition \ref{def:fMV} holds, since it is already satisfied in $(V,u)$.
Hence $[0,u]_V=([0,u], \oplus,\cdot, ^*, \{\alpha\mid \alpha\in [0,1]\},0)$
is an \textit{f}MV-algebra.
\end{enumerate}
\eex

We  list some further examples.

\bex
\begin{enumerate}
\item For any nonempty set $X$, $[0,1]^X$ with point-wise operations of MV-algebra and point-wise products is an \textit{f}MV-algebra.
\item If $(X,\tau)$ is a  topological compact and Hausdorff space  then 
$C(X)=\{f:X\to [0,1]\mid f \mbox{continuous}\}$ with point-wise operations of MV-algebra and point-wise products  
  is an  \textit{f}MV-algebra.
\item  On the Riesz MV-algebra $[0,1]^2$ with  coordinate-wise operations of MV-algebra and coordinate-wise scalar product, we define the internal product by $(x_1,y_1)\cdot (x_2,y_2)=(x_1 x_2, 0)$ for any $x_1,x_2,y_1,y_2\in [0,1]$.  One can  easily prove that the structure obtained in this way is an \textit{f}MV-algebra.
\item    Let   $H$ be an one dimensional Hilbert space over the complex field, $\mathcal{H}$  the set of all bounded Hermitian operators on $H$ and   assume $\mathcal{D}$ is a nonempty subset of $\mathcal{H}$.  Hence $\mathcal{C}^{''}(\mathcal{D})$,  the bicommutant of $\mathcal{D}$  is an \textit{f}-algebra with strong unit $E$, where $E$ is the identity operator.
For more details we refer to \cite[Chapter 8]{RS}.
Then, the unit interval $[\theta,E]$ is an \textit{f}MV-algebra, where $\theta$ is the null operator.
\end{enumerate}
\eex

In \cite{DiND} the proper notion of ideal was defined for PMV-algebras, under the name of \textit{$\cdot$-ideal}. We 
follow a similar approach for \textit{f}MV-algebras.

\bdfn
If $A$ is an \textit{f}MV-algebra then a subset $I\subset A$ is an {\em ideal} if it  is an ideal of the
MV-algebra reduct of $A$. An ideal $I$ of $A$ is a 
 {\em $\cdot$-ideal} if
l if it is a $\cdot$-ideal of the PMV-algebra reduct of $A$.
\edfn

\bfact \label{rem:idPMV}
Assume $A$ is an \textit{f}MV-algebra and $I$ is a $\cdot$-ideal of $A$.
\begin{itemize}
\item[(1)] Since $\alpha x\leq 1x=x$ for any $\alpha\in [0,1]$ and $x\in A$, then
$x\in I$ implies $\alpha x\in I$ for any $\alpha\in [0,1]$.

\item[(2)] When $A$ is unital, $\cdot$-ideals coincide with the ideals of the MV-algebra reduct of $A$, since for any $x \in I$ and any $ y\in A$ we get $x\cdot y \le x\cdot 1=x $.
\end{itemize}
\efact

\noindent For a $\cdot$-ideal   $I$  of an \textit{f}MV-algebra $A$, we define a   congruence relation by:
\vsp{-0.2}
\begin{center}
$x\sim _{I}y$  \quad  if and only if \quad $x\odot y^*\in I$ and $y\odot x^*\in I$
\end{center}
\vsp{-0.2}
By \cite[Remark 2]{LeuRMV}, the quotient of \textit{A} with respect to $\sim_{I}$ is a Riesz MV-algebra.
In order to prove that the quotient is a PMV\textit{f}-algebra, we need to prove that\vsp{-0.2}
\begin{center}
if $x\sim_{I}y$, then $z\cdot x\sim_{I}z\cdot y$ and $x\cdot z\sim_{I}y\cdot z$.
\end{center}
\vsp{-0.2}
But this follows trivially by conditions (f1) and (f2) in Section \ref{subsPMV}.\\
We can now define the product in the quotient algebra $A/\sim_{I}$ by 
\vsp{-0.2}
\begin{center}
$[x]\cdot [y]=[x\cdot y],$
\end{center}
\vsp{-0.2}
and $\quot{A}{I}$ is an   \textit{f}MV-algebra.

We are ready to prove the  subdirect representation theorem with respect to  totally-ordered structures.

\bthm \label{teo:rapprsubdir}
Any \textit{f}MV-algebra \textit{A} is a subdirect product of totally-ordered \textit{f}MV-algebras.
\ethm
\begin{proof}
By \cite[Theorem 5.5]{DiND}, $\mathcal{U}_{\mathbb R}(A)$ is a subdirect product of totally-ordered PMV$f$-algebras,  i.e. there exists a family $\{A_k\}_{k\in K}$ of totally-ordered PMV\textit{f}-algebras and a subdirect representation $\iota: A\to \prod_{k\in K} A_k$. Consequently, there is a family $\{P_k\}_{k\in K}$ of prime $\cdot$-ideals of $\mathcal{U}_{\mathbb R}(A)$ such that $\cap \{P_k|k\in K\}=\{0\}$. By definition, $P_k$ is also a $\cdot$-ideal of $A$, so $A_k$ is a totally-ordered \textit{f}MV-algebras for any $k\in K$.
\end{proof}

\subsection{Categorical equivalence with \textit{fu}-algebras}

\noindent The aim of this section is to prove a categorical equivalence between   \textit{f}MV-algebras and \textit{f}-algebras, generalizing the similar results for  MV-algebras,  PMV-algebra and Riesz MV-algebras. 

 We denote by $\mathbf{fMV}$  the category whose objects are \textit{f}MV-algebras and whose morphisms are homomorphisms of \textit{f}MV-algebras and by  $\mathbf{fuAlg}$ the category whose objects are \textit{fu}-algebras and whose morphisms are  homomorphisms of \textit{f}-algebras that preserve the strong unit.

The functor $\Gamma_f:\mathbf{fuAlg}\to \mathbf{fMV}$ is defined by 
\begin{itemize}
\item $\Gamma_f(V,u)=[0,u]_V$ as in Example \ref{lem:01B},
\item if $g:(V_1,u_1)\to (V_2,u_2)$  is a homomorphism of \textit{fu}-algebras then $\Gamma_f(g)$ is defined as $g|_{[0,u_1]}$.
\end{itemize}

\bthm \label{teo:05}
The functor $\Gamma _f$ establishes a categorical equivalence.
\ethm
\begin{proof}
$\Gamma_f$ is trivially well defined. First, we prove that  for any \textit{f}MV-algebra $A$,
 there exists an fu-algebra (V,u) such that $A\simeq \Gamma_f (V,u)$.

By \cite[Proposition 3]{LeuRMV}  there exists a   Riesz space with strong unit $(V, u)$ such that ${\cal U}_{(\cdot)}(A)\simeq \Gamma_{\mathbb{R}} (V, u)$. Since $A$ has a PMV-algebra reduct, by \cite[Theorems 4.2 and 5.2]{DiND},  $(V,u)$ can be endowed 
with an internal product $\cdot:V\times V\to V$ such that $(V,\cdot)$ is an $f$-ring. Moreover, by easy computation $(V,u)$ satisfies the condition (fa).

One can easily see that the functor $\Gamma_f$ is  faithful: if we have  $h_1,\ h_2: (V_1,u_1)\rightarrow (V_2, u_2)$, such that $\Gamma _f(h_1)=\Gamma _f(h_2)$, then $h_1$ and $h_2$ coincide on the generating set of $(V_1, u_1)$, therefore they coincide on the whole algebra, and $h_1=h_2$.

To end the proof, we show that $\Gamma_f$ is full. Let A, B be \textit{f}MV-algebras such that $A\simeq \Gamma_f(V_1,u_1)$, $B\simeq \Gamma_f(V_2,u_2)$ for $(V_1, u_1)$ and $(V_2, u_2)$ objects  in $\mathbf{fuAlg}$.  If $h:A\rightarrow B$ is a
 homomorphism of \textit{f}MV-algebras then, in particular, $h$ is a  homomorphism of PMV-algebras between the PMV-algebra reducts of $A$ and $B$. By  \cite[Theorem 4.2]{DiND}, there exists a homomorphism of $\ell u$-rings  $h^{\sharp}:(V_1,u_1)\rightarrow (V_2, u_2)$ extending $h$. By \cite[Chapter XV Section 2 Corollary]{Birk}, any morphism of $\ell$-groups between Riesz spaces is linear. It follows that $h^{\sharp}:(V_1,u_1)\rightarrow (V_2, u_2)$ is also linear, so it is a 
homomorphism of \textit{fu}-algebras.
\end{proof}

We are  now able to  complete the diagram:

\begin{center}
\begin{tikzpicture}
  \node (A) {$\mathbf{uR}$};
  \node (B) [below of=A] {$\mathbf{PMV}$};
  \node (C) [right of=A] {$\mathbf{auG}$};
  \node (D) [below of=C] {$\mathbf{MV}$};
  \node (E) [right of=C] {$\mathbf{uRS}$};
  \node (F) [below of=E] {$\mathbf{RMV}$};
\node(H)[left of=A] {$\mathbf{fuAlg}$};
\node(G)[below of=H] {$\mathbf{fMV}$};
\node(J)[right of=E] {$\mathbf{fuAlg}$};
\node(I)[below of=J] {$\mathbf{fMV}$};
 \draw[->] (A) to node [swap] {$\Gamma_{(\cdot)}$} (B);
\draw[->] (H) to node [swap] {$\Gamma_{f}$} (G);
   
\draw[->] (J) to node[swap] {$\mathcal{U}_{(\cdot \ell)}$} (E);
 \draw[->] (A) to node {$\mathcal{U}_{(\cdot \ell)}$} (C);
  \draw[->] (C) to node [swap] {$\Gamma$} (D);
  \draw[->] (B) to node [swap] {$\mathcal{U}_{(\cdot)}$} (D);
\draw[->] (I) to node {$\mathcal{U}_{(\cdot)}$} (F);  
\draw[->] (E) to node [swap] {$\mathcal{U}_{( \ell \mathbb{R})}$} (C);
\draw[->] (H) to node {$\mathcal{U}_{( \ell \mathbb{R})}$} (A);  
\draw[->] (E) to node {$\Gamma_{\mathbb{R}}$} (F);
\draw[->] (G) to node[swap] {$\mathcal{U}_{\mathbb{R}}$} (B);
    \draw[->] (J) to node {$\Gamma_{f}$} (I);
  \draw[->] (F) to node {$\mathcal{U}_{\mathbb{R}}$} (D);
 \node [below=0.5cm, align=flush center] at (D){ Figure 3.};
\end{tikzpicture}
\end{center}

\bcor
Let $A \simeq \Gamma_f(V,u)$. If $A$ is unital then $(V,u)$ is unital and its unit coincide with the strong unit $u$.
\ecor
\begin{proof}
If follows from Lemma \ref{rem:01} and \cite[Proposition 3.1 and Theorem 3.3]{DiND}
\end{proof}
\bcor \label{rem:birk} Any homomorphism of PMV-algebras  between the PMV-algebra reducts of two \textit{f}MV-algebras preserves the scalar product, i.e. it is a morphism of \textit{f}MV-algebras.
\ecor
\begin{proof} By \cite[Chapter XV Section 2 Corollary]{Birk} any morphism of $\ell$-groups between Riesz spaces is linear. Then we apply Theorem \ref{teo:05}. \end{proof}

\bcor \label{cor:prodinfssup}
Let A be an fMV-algebra in which the internal product is commutative. Then $x\cdot y=(x\wedge y)\cdot (y\vee x)$.
\ecor
\begin{proof}
By Theorem \ref{teo:05}, there exists an \textit{fu}-algebra $(V,u)$ such that $A\simeq \Gamma_f(V,u)$. By \cite[Theorem 142.4]{RSZan}, $x\cdot y=(x\wedge y)\cdot (y\vee x)$ for any $x,y \in V$.
\end{proof}

Let \textit{V} be an \textit{f-}algebra. Recall that a  subset \textit{J} of \textit{V} is an $\ell$-ideal if it is a linear subspace that satisfies the following conditions:

($\ell$1) if $x\in V$, $y\in J$, $|x|\le |y|$ then $x\in J$,

($\ell$2)  if $x\in J$, $y\in V$ then $x\cdot y\in J$ and $y\cdot x\in J$.

\bfact \label{teoid:01}
Let (V,u) be an \textit{fu}-algebra, and let  $A=\Gamma(V,u)$. We define the maps
\vsp{-0.2}
\begin{center}
$\Phi: I\mapsto \Phi(I)=\{ x\in V \mid |x|\wedge u\in I\}$ and  $\Psi: H\mapsto \Psi(H) =H \cap [0, u]$
\end{center}
\vsp{-0.2}
for any  $\cdot$-ideal $I$  of $A$ and any  $\ell$-ideal  $H$ of $V$. Then $\Phi$ and $\Psi$ are
 order-isomorphisms from the set of $\cdot$-ideals of $A$ to the set of $\ell$-ideals of $(V,u)$. The proof follows by the similar result for PMV-algebras  and $\ell$-rings  \cite[Theorem 5.1]{DiND} and by the fact that any ideal of $\ell$-ring considered in \cite{DiND} is an $\ell$-ideal by \cite[Chapter XV Section 2 Lemma 1]{Birk}.
\efact

\bdfn

For a $\cdot$-ideal $I$ of  an   \textit{f}MV-algebra $A$, the {\em  \textit{nil}-radical } is
\vsp{-0.2}
\begin{center}
$\sqrt{I}=\{ x\in A \mid x^n\in I \mbox{ for some }n\in \mathbb{N}\},$
\end{center}
\vsp{-0.2}
where $x^n$ is $\overbrace{x\cdot \ldots \cdot x}^{n\,\, times}$.
\edfn
\bprop
If $I$ is a $\cdot$-ideal in a commutative fMV-algebra $A$, then $\sqrt{I}$ is a $\cdot$-ideal in $A$.
\eprop
\begin{proof}
Let $(V,u)$ be the \textit{fu}-algebra such that $A=\Gamma_f(V,u)$. We first remark that $\Psi(\sqrt{\Phi(I)})=\sqrt{I}$, using the notation of Remark \ref{teoid:01}. Indeed:
\vsp{-0.2}
\begin{center}
$x\in \Psi(\sqrt{\Phi(I)}) \quad \Leftrightarrow \quad x\in \sqrt{\Phi(I)}\cap [0,u] \quad \Leftrightarrow$\\
 $ 0\le x\le u \mbox{ and } x\in \sqrt{\Phi(I)} \quad \Leftrightarrow$\\
$ 0\le x\le u \mbox{ and there exists }n\in \mathbb{N} \mbox{ such that } x^n\in \Phi(I) \quad \Leftrightarrow $\\ 
$  0\le x\le u \mbox{ and }|x^n|\wedge u\in I \quad \Leftrightarrow \quad x\in A \mbox{ and } x^n\wedge u=x^n\in I \quad \Leftrightarrow \quad x\in \sqrt{I}.$
\end{center}
\vsp{-0.2}
By Remark \ref{teoid:01}, $\Phi(I)$ is an $\ell$-ideal, then $\sqrt{\Phi(I)}$ is an $\ell$-ideal by
 \cite[Proposition 4.2]{HuDePID} and again by Remark \ref{teoid:01} $\Psi(\sqrt{\Phi(I)})$ is a $\cdot$-ideal.
\end{proof}

\section{Classes of \textit{f}MV-algebras}\label{sec2}
\noindent  In the theory of \textit{f}-algebras,  the Archimedean, the semisimple, the semiprime and the formally real structures are proper subclasses studied in the literature. We introduce and briefly investigate the same notions  for \textit{f}MV-algebras.
Our main goal is achieved in Section \ref{frplus}, where we characterize the quasi-variety of \textit{f}MV-algebras generated by $[0,1]$.

\subsection{Semisimple and Archimedean \textit{f}MV-algebras}
\noindent It is known that  Archimedean MV-algebras and semisimple MV-algebras coincide. The same holds for Riesz MV-algebras, since a Riesz MV-algebra has the same congruences as its MV-algebra reduct.  The theory of \textit{f}MV-algebras has this property no longer.

We recall that an \textit{f}-algebra is {\em Archimedean} if $na\le b$ for any $n\in \mathbb{N}$ implies $a\le 0$. The notion of \textit{semisimple} \textit{f}-algebra is present in literature in different forms, due to the different kind of radicals that can be defined on a \textit{f}-ring (see for example \cite[Section 8.6]{BigKim}). We will follow \cite{MW}, and define an \textit{f}-algebra $\ell$-semisimple if the intersection of all maximal $\ell$-ideals is $\{ 0 \}$.

\bdfn
Let $A$ be an \textit{f}MV-algebra and let $(V,u)$ be an \textit{fu}-algebra such that $A\simeq\Gamma_f(V,u)$.
We say that the  \textit{f}MV-algebra $A$ is  \textit{Archimedean} ($\ell$-\textit{semisimple}) if the \textit{fu}-algebra 
$V$ is Archimedean ($\ell$-semisimple). 
\edfn

 We note that an  \textit{f}MV-algebra $A$ is Archimedean if and only if its MV-algebra reduct $\mathcal{U}_{( \mathbb{R},\cdot )}(A)$ is a semsimple MV-algebra. 

\bprop
Let $A$ be an Archimedean \textit{f}MV-algebra. Then $A$ is commutative with respect to the internal product.
\eprop
\begin{proof}
Let $(V,u)$ be the \textit{fu}-algebra such that $A\simeq \Gamma_f(V,u)$. By definition, $A$ is Archimedean if and only if $(V,u)$ is Archimedean. Therefore, the result follows from \cite[Theorem 140.10]{RSZan}.
\end{proof}

\blem Any $\ell$-semisimple \textit{f}MV-algebra is Archimedean.
\elem
\begin{proof}
 If we denote by $\ell Rad(A)$ the intersection of all maximal $\cdot$-ideal for the \textit{f}MV-algebra $A$, since any $\cdot$-ideal is an ideal for the reduct $\mathcal{U}_{(\mathbb{R}, \cdot)}(A)$ it follows that $Rad(A) \subseteq \ell Rad(A)$. Then, if $A$ is $\ell$-semisimple, so it is  $\mathcal{U}_{(\mathbb{R}, \cdot)}(A)$, and $A$ is Archimedean.
\end{proof}

\bfact
If $A$ is a unital \textit{f}MV-algebra then $A$ is Archimedean if and only if $A$ is $\ell$-semisimple. This result is a direct consequence of  Remark \ref{rem:idPMV} (2):  in this case any ideal is a $\cdot$-ideal.
 \efact 

 For unital \textit{f}MV-algebras, the two classes coincide. Our main interest focuses on this case, but the general relation between
Archimedeanity and semisimplicity in the non-unital case will make the subject of future studies.

\bprop \label{pro:pmv}
Let $A_1$, $A_2$ be semisimple and unital   \textit{f}MV-algebras, and let $f:A_1 \rightarrow A_2$ be a homomorphism of Riesz MV-algebras. Then $f$ is also a homomorphism of PMV-algebras.
\eprop
\begin{proof}
By Theorem \ref{teo:05} there exist two unital and Archimedean \textit{f}-algebras $(V_1, u_1)$ and $(V_2,u_2)$ such that $A=\Gamma _f(V_1,u_1)$ and $B=\Gamma _f(V_2,u_2)$. The map $f^{\sharp}$ that extend $f$ is a positive map between Archimedean and unital $f$-algebras such that $f^{\sharp}(v)=u$. Therefore by \cite[Corollary 5.5]{HuDeP},  $f^{\sharp}$ is also an algebra homomorphism, and $f:A_1 \rightarrow  A_2$ is a homomorphism of PMV-algebras, as well as a homomorphism of Riesz MV-algebras.
\end{proof}

\subsection{Semiprime   \textit{f}MV-algebras}

\noindent We recall that a {\em semiprime}  \textit{f}-algebra  is an algebra without nilpotents \cite[Chapter 142]{RSZan}. 

\bdfn
An element $x$ of an \textit{f}MV-algebra (\textit{f}-algebra) is {\em nilpotent} if there exists an positive integer $n$ such that $x^n=\underbrace{x\cdot \ldots \cdot x}_{n\,\, times}=0$.
 An \textit{f}MV-algebra $A$  (\textit{f}-algebra $V$) is called \em{semiprime} if $x^n=0$ implies $x=0$.
\edfn

\bfact 
One can easily see that an \textit{f}MV-algebra $A$ (\textit{f}-algebra $V$)  is semiprime if and only if $x\cdot x=0$ implies $x=0$ for any $x\in A$ ($x\in V$). 
\efact

\bprop \label{pro:01s}
 Assume $A$ is an \textit{f}MV-algebra and $(V,u)$  is an \textit{fu}-algebra such that $A\simeq\Gamma_f(V,u)$.  
Then $A$  is semiprime if and only if  $(V,u)$ is semiprime.
\eprop
\begin{proof}
In one direction it is obvious. For the other direction, in the case of positive elements the proof is a simple matter of computation; then the result follows from the positive element case by \cite[Theorem 142.1(iv)]{RSZan}, that is the remark that $x^+\cdot x^-=0=x^-\cdot x^+$.
\end{proof}

\bfact
Since the property of being semiprime depends only on the underlying \textit{f}-ring structure, a  similar  result  holds for PMV\textit{f}-algebras and \textit{f}-rings. We state it without proof.
\efact

\bprop
 Assume $P$ is a commutative and unital PMV-algebra and $(R,u)$  is an \textit{fu}-ring such that $P\simeq\Gamma_{(\cdot)}(R,u)$.  
Then $P$  is a PMV$^+$-algebra  if and only if  $(R,u)$ is semiprime.
\eprop

\bcor
Let $A$ be a semisimple and unital \textit{f}MV-algebra. Then $A$ is semiprime.
\ecor
\begin{proof}
Let $(V,u)$ be the \textit{fu}-algebra such that $A=\Gamma_f(V,u)$. By hypothesis, $A$ is Archimedean, so $V$ is also Archimedean and $u$ is a unit for the product. By \cite[Theorem 142.5(iii)]{RSZan} it follows that  $(V,u)$ is semiprime, so  $A$  is semiprime by Proposition \ref{pro:01s}.
\end{proof}

\bcor
Let $A$  be a semiprime \textit{f}MV-algebra. Then for any $x,y\in A$, $x^2=y^2$ if and only if $x=y$.
\ecor
\begin{proof}
The result follows directly by Proposition \ref{pro:01s} and \cite[Theorem 142.3(ii)]{RSZan}.
\end{proof}

We shall further analyze  the commutative, unital and semiprime algebras in the following sections.

\subsection{Formally real \textit{f}MV-algebras}\label{formreal}

\noindent While the varieties of Abelian $\ell$-groups and Riesz spaces are generated by $\mathbb{R}$ endowed with the corresponding structure, this is not true for the varieties of \textit{f}-rings and \textit{f}-algebras and counterexamples can be found in \cite{HI, Madden}.  Consequently,  algebras from HSP$([0,1])$ were called {\em formally real}. 

The same facts remain true when applying the $\Gamma$ functors (see Figure 3). The varieties of MV-algebras and Riesz MV-algebras are generated by the real interval $[0,1]$ endowed with the corresponding structure (see \cite{CDM} and, respectively \cite{LeuRMV}), while the varieties of (commutative and unital) PMV-algebras  and \textit{f}MV-algebras are not.
For PMV-algebras, a counterexample was given in \cite{HC}. This example also stands for    \textit{f}MV-algebras and we 
briefly  recall it in the next example.

\bex
In  \cite[Example 3.14]{HC} the authors define a totally-ordered finite monoid $S$ and they consider 
the set $ F[S]= \{ r_1 X^{s_1} +\ldots +r_nX^{s_n} \mid n \in \mathbb{N},\ r_i\in F\ s_i\in S\}$, where $F$   is the ordered field of real numbers.  They further identify $X^\top$ with $0$, where $\top$ is the greatest element of $S$ and they denote 
$F[S]_h$ the quotient obtained in this way, which is an \textit{f}-ring.  Hence the interval  $[\mathbf{0},\mathbf{1}]$ of $F[S]_h$
is a PMV-algebra  that does not satisfy the  following  identity:

$(x_1 \cdot z_1 \ominus y_1 \cdot z_2) \wedge (x_2 \cdot z_2 \ominus y_2 \cdot z_1)\wedge (y_1 \cdot y_2 \ominus x_1 \cdot x_2)=0$.\\
Since the identity holds in the real interval $[0,1]$, it provides the intended count-example in the context of PMV-algebras. 
We only have to note that  $F[S]_h$ is in fact an $\textit{f}$-algebra so $[\mathbf{0},\mathbf{1}]$ is an \textit{f}MV-algebra that does not belong to HSP$([0,1])$.
\eex

\bdfn
Following \cite{HI}, we will call an \textit{f}MV-algebra (PMV-algebra) \textit{formally real} if it belongs to HSP$([0,1])$. We denote by $\mathbb{FR}$ the class of formally real   \textit{f}MV-algebras.
\edfn

\bfact
In general, a formally real $f$-ring is not unital. For PMV-algebras and   \textit{f}MV-algebras the situation is different: by \cite{DiND} the unit is the greatest element of the algebras and it belongs to the language. Then the condition for the unit is not existential but universal, that is $1\cdot x=x\cdot 1=x$ for any $x$.
\efact

By well-known results of universal algebra (see for example \cite{Gr}), the free \textit{f}MV-algebra in  ${\mathbb F}{\mathbb R}$ exists and its elements are term functions defined on $[0,1]$.  More precisely,  the language of \textit{f}MV-algebras is 

${\cal L}_f=\{\oplus, \cdot, ^*, 0,\}\cup\{\delta_{\alpha}\mid \alpha \in [0,1]\}$,\\
where $\delta_{\alpha}$ is a unary operation that is 
interpreted by $x\mapsto \alpha x$ for any $\alpha \in [0,1]$. For any $n\geq 1$, let $X=\{x_1,\ldots, x_n\}$ and  assume
$Term_n$ is the set of ${\cal L}_f$-terms with variables from $X$. We  
   denote by $FR_n$ the free \textit{f}MV-algebra in ${\mathbb F}{\mathbb R}$ with $n$ free generators. It follows that
\vsp{-0.1}
\begin{center}
$FR_n=\{\widetilde{t}\mid t\in Term_n, \,\, \widetilde{t}: [0,1]^n\to [0,1] \mbox{ is the term function of } t\}$.
\end{center}
\vsp{-0.1}
In order to characterize $FR_n$ we give the following definition. 

\bdfn\label{pwpdef}
A {\em piecewise polynomial function} defined on the $n$-cube
is a continuous function   $f:[0,1]^n\rightarrow [0,1]$  such that there  exists a finite number of polynomials  $f_1$, $\ldots$, $f_k\in {\mathbb R}[x_1,\ldots, x_n]$  with the property that  $f(a_1,\ldots, a_n)=f_i(a_1,\ldots, a_n)$ 
  for any $(a_1,\ldots, a_n)\in  [0,1]^n$  and for some $i\in \{1,\ldots, k\}$.

The polynomials $f_i$ are called {\em the components} of $f$.
\edfn

\bprop\label{ppf}
The elements of $FR_n$ are piecewise polynomial functions defined on the $n$-cube.
\eprop
\begin{proof}
Let $t$ be a term in $Term_n$. The result will be proved by structural induction on $t$. We recall that any terms in $Term_n$ is built from the language $\mathcal{L}_f$ and the set of variables $\{ x_i \}_{i\le n}$.

 If $t=x_i$ for some $i\le n$, then $\widetilde{t}=\pi_i^n$, the $i^{th}$ projection and it is trivially a piecewise polynomial function.

 If $t=t_1^*$, then $\widetilde{t}=\widetilde{t_1}^*$. By induction hypothesis there exists an integer $h$ and some polynomials $q_1, \ldots , q_h \in \mathbb{R}[x_1, \ldots ,x_n]$ such that for any point in the $n$-cube, $\widetilde{t_1}$ coincides with one of them. Then $1-q_1, \ldots, 1-q_h$ are the components of $\widetilde{t}$.

 If $t= t_1 \oplus t_2$, let $q_{1}$,$\ldots$, $q_{m}$  be the  
components of $\widetilde{t_1}$ and 
$p_{1}$,$\ldots$,$p_{k}$ be the components of $\widetilde{t_2}$. Then $\widetilde{t}$ is defined by the set of polynomials
$\{1\}\cup \{s_{ij}\}_{i,j}$, where $s_{ij}=1-q_i+p_j$ for any $i\in\{1,\ldots ,m\}$ and $j\in\{1,\ldots ,k\}$.

 If $t=\delta_{\alpha}(t_1)$ for some $\alpha \in [0,1]$ and  $q_{1}$,$\ldots$,$q_{s}$  are  the components of $\widetilde{t_1}$, then $\alpha q_{1}$,$\ldots$,$\alpha q_{s}$ are the components of $\widetilde{t}$.

If $t= t_1 \cdot t_2$, let $q_{1}$,$\ldots$, $q_{m}$  be the  
components of $\widetilde{t_1}$ and 
$p_{1}$,$\ldots$,$p_{k}$ be the components of $\widetilde{t_2}$. Then $\widetilde{t}$ is defined by the polynomials $q_i\cdot p_j$, for any $i\in\{1,\ldots ,m\}$ and $j\in\{1,\ldots ,k\}$.
\end{proof}

The converse of the above proposition is related to the Pierce-Birkhoff conjecture \cite{BP,HI} and it will be analyzed in Section \ref{pbc}.

\subsection{FR$^+$-algebras}\label{frplus}

\noindent The class of FR$^+$-algebras is the quasi-variety of \textit{f}MV-algebras generated by $[0,1]$. In order to characterize this class, we follow  the ideas from \cite{Mon+}, where a similar investigation is done for PMV-algebras.

The  FR$^+$-algebras  are the core of our  development, from the algebraic point of view. We prove a representation theorem w.r.t. totally-ordered structures, as well as a representation theorem via ultrapowers of $[0,1]$. In Section \ref{pblogic} we develop a logical system that has FR$^+$-algebras as models and in Section \ref{pbc} we connect them this logic with the Pierce-Birkhoff conjecture.

We recall that a PMV-algebra $A$ is a   $PMV^+$-algebra  \cite{Mon+} if it is  unital, commutative and it satisfies the condition:

(sp) $x\cdot x=0$ implies $x=0$ for any $x\in A$.

\bthm\label{teo:qvsemiprime}
 For an \textit{f}MV-algebra $A$ the following are equivalent:\\
(1) $A\in $ ISP$([0,1])$, \\
(2) $A$ is a  unital, commutative and semiprime fMV-algebra,\\
(3) $\mathcal{U}_{\mathbb{R}}(A)$, the PMV-algebra reduct of $A$, is a PMV$^+$-algebra,\\
(4)  $A$ is in $\mathbb{FR}$ and   $x\cdot x=0$ implies $x=0$ for any $x\in A$.
\ethm
\begin{proof}
$(1)\Leftrightarrow (2)$ One direction is trivial, since any element in ISP$([0,1])$ is unital, commutative and semiprime. For the other direction, let $\mathbf{V}$ be the class of unital, commutative and semiprime   \textit{f}MV-algebras. If $A\in \mathbf{V}$, then $\mathcal{U}_{\mathbb{R}}(A)$ is a $PMV^+$-algebras, and by \cite[Corollary 4.4]{Mon+} it belongs to ISP$([0,1])$, therefore $A$ is a subalgebra of a direct product of copies of $[0,1]$, and the direct product is trivially a \textit{f}MV-algebra. The map that gives the inclusion is a homomorphism of PMV-algebras, therefore it is a homomorphism of  \textit{f}MV-algebras, and $A$  belongs to ISP$([0,1])$ as    \textit{f}MV-algebra.\\
$(2)\Leftrightarrow (3)$ is obvious.\\
$(1)\Rightarrow (4)$  Since ISP$([0,1])\subseteq$ HSP$([0,1])$, it is straightforward.\\
$(4) \Rightarrow (2)$ Any formally real \textit{f}MV-algebra is unital and commutative, and the quasi-identity in the hypothesis characterizes semiprime \textit{f}MV-algebras.
\end{proof}

\bdfn
Following \cite{Mon+} we denote  $\mathbb{FR}^+=$ISP$([0,1])$, the quasi-variety generated by the standard \textit{f}MV-algebra $[0,1]$.  An \textit{f}MV-algebra $A$ is a FR$^+$-algebra if it is unital, commutative and semiprime.
\edfn 

\bprop\label{free}
$FR_n$ belongs to $\mathbb{FR}^+$, for any positive integer $n$.
\eprop
\begin{proof}
By Theorem \ref{teo:qvsemiprime} we just need to prove that $FR_n$ is semiprime. Let $g\in FR_n$ such that $g^2=\mathbf{0}$, the null function. This means that $g(\mathbf{x})\cdot g(\mathbf{x})=0$ for every $\mathbf{x}\in [0,1]^n$. Since $g(\mathbf{x})\in [0,1]$, the above condition implies $g(\mathbf{x})=0$ for any $\mathbf{x}\in [0,1]^n$, and $g=\mathbf{0}$.
\end{proof}

For FR$^+$-algebras we prove the subdirect representation theorem w.r.t. totally-ordered FR$^+$-algebras, as well as the representation as an algebra of $^*[0,1]$-valued functions.

\bthm \label{pro:+sp}
Any FR$^+$-algebra is a subdirect product of totally-ordered FR$^+$-algebras.
\ethm
\begin{proof}
It is similar to the proof of Theorem \ref{teo:rapprsubdir}, by \cite[Theorem 3.9]{HC}.
\end{proof}

\bthm
For any FR$^+$-algebra $A$ there exists an ultrapower $\ ^{*}[0,1]$ of $[0,1]$ and a set $I$  such that $A$  is embedded  in $(\ ^{*}[0,1])^I$.
\ethm
\begin{proof}
We first prove the result for totally-ordered FR$^+$-algebras.  By hypothesis, the PMV-algebra reduct of $A$ is a $PMV^+$-chain, therefore by \cite[Corollary 4.4]{Mon+} $A$ can be embedded, as a PMV-algebra, in an ultrapower of $[0,1]$. Since any PMV-algebras homomorphism is also a Riesz MV-algebras homomorphism, the embedding is between \textit{f}MV-algebras.\\
By Proposition \ref{pro:+sp}, any FR$^+$-algebra $A$ is subdirect product of a family $\{A_i\}_i$ of totally-ordered FR$^+$-algebras. We know that any $A_i$ is embedded in $U_i$, an ultrapower of $[0,1]$. By \cite[Proposition 3.1.4 and Corollary 4.3.13]{CK}  there exists an ultrapower $U$ of $[0,1]$ such that every $U_i$ is embedded in  $U$. It follows that  $A$ is embedded in $\Pi_i U$.
\end{proof}

\section{Logic for unital and commutative   \textit{f}MV-algebras}\label{sec3}

\noindent In this section, we define the propositional calculus  $\mathcal{FMVL}$, that has unital and commutative \textit{f}MV-algebras as models, as well as the propositional   calculus  $\mathcal{FMVL^+}$ that has FR$^+$-algebras as models. For the latter we prove standard completeness w.r.t. $[0,1]$.

\subsection{The propositional calculus  $\mathcal{FMVL}$}
\noindent The language of the propositional logic $\mathcal{FMVL}$ consists  of:\\
(i) a countable set of propositional variables $v_1, v_2, \ldots$; \\
(ii) the binary  connectives  are $\rightarrow$ and  $\cdot$;\\
(iii) the unary connective $\neg$;\\
(iv)  a family of unary connectives  $\{\nabla_\alpha\mid \alpha \in [0,1]\}$;\\
(v)  the parentheses $($ and $)$.\\
 Formulas, theorems, deductions, proof are defined as usual.
The axioms for the logic of unital and commutative \textit{f}MV-algebras will be the following:\\
(L1) $\varphi \rightarrow (\psi \rightarrow \varphi)$ \\
(L2) $(\varphi \rightarrow \psi)\rightarrow((\psi \rightarrow \chi)\rightarrow(\varphi \rightarrow \chi))$ \\
(L3) $(\varphi \vee \psi)\rightarrow (\psi \vee \varphi)$ \\
(L4) $(\neg \psi \rightarrow \neg \varphi)\rightarrow (\varphi \rightarrow \psi)$ \\
(R1) $\nabla_{\alpha}(\varphi \rightarrow \psi)\leftrightarrow (\nabla_{\alpha}\varphi \rightarrow \nabla_{\alpha}\psi)$ \\
(R2) $\nabla_{(\alpha \odot \beta ^*)}\varphi \leftrightarrow (\nabla_{\beta}\varphi \rightarrow \nabla_{\alpha}\varphi)$ \\
(R3) $\nabla_{\alpha}(\nabla_{\beta}\varphi)\leftrightarrow \nabla_{\alpha \cdot \beta} \varphi $ \\
(R4) $\nabla_1 \varphi \leftrightarrow \varphi $ \\
(P1) $(\chi \cdot (\varphi \ominus \psi))\leftrightarrow ((\chi \cdot \varphi)\ominus (\chi \cdot \psi)) $ \\
(P2) $(\varphi \cdot (\psi \cdot \chi))\leftrightarrow ((\varphi \cdot \psi)\cdot \chi) $ \\
(P3) $\varphi \rightarrow (\varphi \cdot (\varphi\rightarrow\varphi)) $ \\
(P4) $\varphi \cdot \psi \rightarrow \varphi  $ \\
(P5) $\varphi \cdot \psi \leftrightarrow \psi \cdot \varphi $ \\
(A1) $ \Delta_{\alpha}(\varphi \cdot \psi ) \leftrightarrow (\Delta_{\alpha }\varphi\cdot \psi) $ \\
(A2) $ \Delta_{\alpha}(\varphi \cdot \psi)\leftrightarrow (\varphi \cdot \Delta_{\alpha}\psi),$\\
where  $\varphi \ominus \psi$ means $\neg (\varphi \rightarrow \psi)$, $\Delta_{\alpha} \varphi$ means $\neg \nabla_{\alpha}(\neg \varphi)$  and the only deduction rule will be the \textit{Modus Ponens}. The set of formulas in $\mathcal{FMVL}$ will be denoted by $F_{fMV}$.

Note that (L1)-(L4) are the axioms of \L ukasiewicz logic.

For a subset $\Theta \subseteq F_{fMV}$ we define an   equivalence relation as follows:
\vsp{-0.2}
\begin{center}
$\varphi \equiv_{\Theta} \psi \mbox{\quad iff \quad} \Theta \vdash \varphi \rightarrow \psi \mbox{ and }\Theta \vdash \psi \rightarrow \varphi.$
\end{center}
\vsp{-0.2}
We define the following operation on $F_{fMV}/\equiv_{\Theta}$:

$\bullet$ $[\varphi ]_{\Theta}^*=[\neg \varphi]_{\Theta}$, \quad $[\varphi]_{\Theta} \rightarrow [\psi]_{\Theta}= [\varphi \rightarrow \psi]_{\Theta}$;

$\bullet$ $[\varphi]_{\Theta} \oplus [\psi]_{\Theta}=[\neg \varphi\rightarrow \psi]_{\Theta}$, \quad $[\varphi]_{\Theta}\cdot [\psi]_{\Theta}=[\varphi \cdot \psi]_{\Theta}$;

$\bullet$ $\alpha [\varphi]_{\Theta}=[\Delta_{\alpha}\varphi]_{\Theta}$

$\bullet$ $1_{\Theta}=T_{\mathcal{FMVL}}(\Theta)$, \quad $0_{\Theta}=1^*_{\Theta}$.

\bthm
The structure $fMVL(\Theta)=((F_{fMV}/\equiv_{\Theta}, \oplus, *, \cdot, 0_{\Theta}), \Phi)$, where $\Phi:[0,1]\times F_{fMV} \rightarrow F_{fMV},\ \Phi(\alpha, [\varphi])=[\Delta_{\alpha} \varphi]$,  is a unital and commutative fMV-algebra.
\ethm
\begin{proof}
By \cite[Proposition 5]{LeuRMV}, this structure is a \textit{Riesz} MV-algebra. By \cite[Lemma 3.21]{HC}  it is a unital and commutative PMV-algebra, and by axioms (A1) and (A2) we get the algebra's relations between internal product and scalar product. Trivially any unital PMV-algebra is an PMV\textit{f}-algebra. 
\end{proof}

In the sequel, we will denote by $fMVL$ the Lindenbaum-Tarski algebra of the logic, that is $fMVL(\Theta)$ with $\Theta=\emptyset$.

An appropriate semantic for the logic has unital and commutative \textit{f}MV-algebras as models. Let $A$  be an   \textit{f}MV-algebra; an evaluation is a function $e: F_{fMV}\rightarrow A$ such that:

(e1) $e(\varphi \rightarrow \psi)= e(\varphi)^* \oplus e(\psi)$,

(e2) $e(\neg \varphi =e(\varphi)^*)$,

(e3) $e(\nabla_{\alpha}\varphi)=(\alpha e(\varphi)^*)^*$,

(e4) $e(\varphi \cdot \psi)=e(\varphi)\cdot e(\psi)$,\\
for any formulas $\varphi$ and $\psi$ of $\mathcal{FMVL}$.    

We have the following completeness theorem.

\bthm
Let $\Theta$ be a set of formulas and $\varphi$ a formula in $F_{fMV}$. The following are equivalent:\\
(1) $\Theta \vdash \varphi$,\\
(2) $\Theta \models _A \varphi$ for any \textit {fMV}-algebra A,\\
(3) $\Theta \models _A \varphi$ for any totally-ordered \textit {fMV}-algebra A,\\
(4) $[\varphi]_{\Theta}=1_{\Theta}$ in $fMVL(\Theta)$.
\ethm
\begin{proof}
(2) $\Leftrightarrow$ (3)  follows by  Theorem \ref{teo:rapprsubdir}. The rest of the proof  is straightforward.
\end{proof}

\bprop
$\mathcal{FMVL}$ is a conservative extension of $\mathcal{L}$, the infinite valued \L ukasiewicz logic.
\eprop
\begin{proof}
Let $\varphi$ be a formula in \L ukasiewicz logic which is a theorem in $\mathcal{FMVL}$. Then by completeness, $\varphi $ is a tautology in the standard \textit{f}MV-algebra $[0,1]$, and since $\varphi$ does not involves $\cdot$ and $\nabla _{\alpha}$, it is a tautology in the standard MV-algebra $[0,1]$. Then by the completeness of \L ukasiewicz logic, $\varphi$ is a theorem in $\mathcal{L}$.
\end{proof}

As for \L ukasiewicz logic,  the deduction theorem holds in its local form.

\bfact
Let $\Theta$ be a non-empty subset of $F_{fMV}$, and $\varphi \in F_{fMV}$. For any $\psi \in F_{fMV}$ we have
$$\Theta \cup \{ \varphi \} \vdash \psi \mbox{\quad iff \quad} \Theta \vdash\varphi\rightarrow \cdots(\varphi {\rightarrow} \psi),$$
where $\varphi$ appears $n$ times for some $n\geq 1$.
The proof is similar  with the one for \L ukasiewicz logic \cite[Proposition 4.6.4]{CDM}.
\efact

\subsection{The propositional calculus  $\mathcal{FMVL}^+$}\label{pblogic}

\noindent  The propositional calculus   $\mathcal{FMVL}^{+}$ is obtained  from  $\mathcal{FMVL}$  by adding the deduction rule
\begin{center}
{\em Semiprime}: $ \displaystyle{\frac{\neg (\varphi \cdot \varphi)}{\neg \varphi}}.$
\end{center}

\bthm \label{coplus}
Let $\Theta$ be a set of formulas in $\mathcal{FMVL}^{+}$ and $\varphi $ a formula. Then the following are equivalent:\\
(1) $\Theta\vdash \varphi$,\\
(2) $e(\varphi)=1$ for any  $e$ $[0,1]$-model of $\Theta$,\\
(3) $e(\varphi)=1_A$ for any algebra  $A$ from  $\mathbb{FR}^+$  and for any $e$ $A$-model of $\Theta$.
\ethm
\begin{proof}
$(1)\Rightarrow (2)$ It is straightforward, since {\em Modus Ponens} and {\em Semiprime}  preserve  tautologies.\\
$(2)\Rightarrow (3)$ It follows by Theorem \ref{teo:qvsemiprime} \\
$(3)\Rightarrow (1)$ Since (3) holds for the algebras from $\mathbb{FR}^+$, then (1) follows directly by the definition of $1_{\Theta}.$
\end{proof}

It follows  that the models of  $\mathcal{FMVL}^{+}$ are $FR^+$-algebras. Note that 
our system is an extension of the system $PL^\prime$,  defined in \cite{HC}, that has PMV$^+$-algebras as models.

\bprop
$\mathcal{FMVL}^{+}$ is a conservative extension of $\mathcal{L}$.
\eprop
\begin{proof}
It is similar to the proof for $\mathcal{FMVL}$.
\end{proof}

\bfact
$\mathcal{FMVL}^{+}$ does not satisfy the deduction theorem in the same form of $\mathcal{FMVL}$: the counterexample in  \cite[Corollary 3.19]{HC} applies. 
\efact

 The appropriate  semantic for $\mathcal{FMVL}^{+}$ is the one with unital, commutative and semiprime   \textit{f}MV-algebras, that is, algebras in $\mathbb{FR}^+$, as models with the usual definition for evaluations. 

\bfact
Let $n\geq 1$ be a natural number. 
The Linbenbaum-Tarski algebra with $n$ variables, denoted  by $fMVL^{+}_n$,   is isomorphic with   the free FR$^+$-algebra with $n$ free variables. By Propositions \ref{ppf} and \ref{free}, we infer that  $fMVL^{+}_n\simeq FR_n$ and we know that  the elements of $FR_n$ are $[0,1]$-valued piecewise polynomial functions defined on $[0,1]^n$.   
\efact

In the following section we shall investigate the answer to the following question:
{\em does $FR_n$ contain \underline{all} the  $[0,1]$-valued piecewise polynomial functions defined on $[0,1]^n$}?

The following result is a preliminary step.

\bprop\label{formfct}
If $p:[0,1]^n\to {\mathbb R}$  is  a polynomial with  real  coefficients then  there exists a formula $\varphi$ with $n$ variables of  $\mathcal{FMVL}^{+}$  such that  $((p\vee {\mathbf 0})\wedge {\mathbf 1})$ coincides with the term function associated to $\varphi$. 
\eprop
\begin{proof}
  The  proof is similar with  the one of  \cite[Proposition 7.6]{LeuRMV}.
 In the following we denote the function $p \mapsto ((p\vee {\mathbf 0})\wedge {\mathbf 1})$ by $\varrho$, and by $\widetilde{\varphi}$ the term function associated to the formula $\varphi$.
Let $p:[0,1]^n\to {\mathbb R}$  be a polynomial function. Let $k$ be the degree of $p$. Then it follows
$$ p(x_1, \ldots , x_n)= \sum_{i_1+\ldots +i_n\le k}c_{i_1,\ldots ,i_n }x_1^{i_1}\cdots x_n^{i_n},$$
where $c_{i_1,\ldots ,i_n } \in \mathbb{R}$ for any choice of the indexes. We notice that any $c_{i_1,\ldots ,i_n }$ can be written as a sum of a finite number of elements in $[-1,1]$, then we  assume that
$$ p(x_1,\ldots, x_n)=r_my_m+\cdots r_{p+1}y_{p+1}+r_{p}+\cdots+r_1 $$
where $m\geq 1$ and $p\geq 0$ are natural numbers, $p\leq m $,  $r_j\in [-1,1]\setminus\{0\}$ for any $j\in\{1,\ldots,\ m\}$ and
 $y_j\in\{x_1^{i_1}\cdots x_n^{i_n} \mid i_1+\ldots +i_n \le k\}$ for any $j\in\{p+1,\ldots,\ m\}$.\\
\noindent We prove the theorem by induction on $m\geq 1$.
In the sequel we denote by $\textbf{x}$ an element $(x_1,\cdots, x_n)$ from  $ [0,1]^n$.\\
\noindent {\em Initial step} $m=1$. We have
$p(\textbf{x})=r$ for any $\textbf{x}\in [0,1]^n$ or  $p(\textbf{x})=rx_1^{i_1}\cdots x_n^{i_n}$ for any $\textbf{x}\in [0,1]^n$ where $r\in  [-1,1]\setminus\{0\}$ and $\{ i_1, \ldots , i_n \}$ is a suitable set of index.

If  $r\in[-1,0)$ then $\varrho ( p)=0$  so $\varrho\circ p=\widetilde{\varphi}$ for $\varphi=\overline{0}$.

If $r\in (0,1]$ then  $p=\varrho\circ p$.  It follows that $p=\widetilde{\varphi}$ where $\varphi=\delta_r(\overline{0}^*)$
if $p(\textbf{x})=r$ for any $\textbf{x}\in [0,1]^n$, and  $\varphi=\delta_r (v_1^{i_1}\cdot \ldots \cdot v_n^{i_n})$
if   $p(\textbf{x})=rx_1^{i_1}\cdots x_n^{i_n}$  for any $\textbf{x}\in [0,1]^n$.\\
\noindent{\em Induction step}. We take  $p=g+h$ where $\varrho\circ g=\widetilde{\varphi_1}$ for some formula $\varphi_1$ and there is
 $r\in  [-1,1]\setminus\{0\}$ and a suitable choice of index for $i_1, \ldots , i_n$ such that
$h(\textbf{x})=r$  for any $\textbf{x}\in [0,1]^n$, or  $h(\textbf{x})=rx_1^{i_1}\cdots x_n^{i_n}$ for any $\textbf{x}\in [0,1]^n$.
We consider two cases.\\
\noindent{\em Case 1.} If   $r\in (0,1]$ then $h:[0,1]^n\to [0,1]$ so
$\varrho\circ p =((\varrho\circ g)\oplus h)\odot (\varrho\circ (1+g))$ by \cite[Lemma 10]{LeuRMV}.
  Following the initial step, there is a formula $\varphi_2$ such that $h=\widetilde{\varphi _2}$. We notice that  $1+g=1-(-g)$ and the induction hypothesis holds for $(-g)$, then
there is a  formula $\varphi _3$ such that $\varrho\circ (-g)=\widetilde{\varphi _3}$. It follows by \cite[Lemma 10]{LeuRMV},
$\varrho\circ (1+g)=1-\widetilde{\varphi _3}=\widetilde{\varphi _3}^*$.
We get $\varrho \circ p=\widetilde{\varphi}$ where $\varphi=(\varphi _1 \oplus  \varphi _2)\odot \varphi _3^*$.\\
\noindent{\em Case 2.} If $r\in [-1,0)$, then   $g+h=(g-1)+(1+h)$ and $1+h:[0,1]^n\to [0,1]$. By  \cite[Lemma 10]{LeuRMV} we get
\vsp{-0.2}
\begin{center}
$\varrho\circ p =((\varrho\circ (g-1))\oplus (1+h))\odot (\varrho\circ g).$
\end{center}
\vsp{-0.2}
\noindent Following the initial step, there is a formula $\varphi _2$ such that $-h=\widetilde{\varphi _2}$, so $1+h=1-(-h)=\widetilde{\varphi_2}^*$.
In the sequel we have to find a formula $\varphi_3$ that corresponds to $\varrho\circ (g-1)$, where
\vsp{-0.2}
\begin{center}
$ g(\textbf{x})=r_my_m+\cdots +r_{p+1}y_{p+1}+r_{p}+\cdots+r_1 $
\end{center}
\vsp{-0.2}
with  $r_j\in [-1,1]\setminus\{0\}$ for any $j\in\{1,\ldots, m\}$ and  $y_j$ in $\{ x_1^{i_1}\cdots x_n^{i_n} \mid i_1+\ldots + i_n \le k\}$ for any $j\in\{p+1,\ldots, m\}$.\\
\noindent{\em Case 2.1.} If  $r_j\leq 0$ for any $j\in\{1,\ldots, m\}$ then $g-1\leq 0$, so $\varrho\circ (g-1)=0= \widetilde{\varphi_3}$ with
 $\varphi_3=0$.\\
\noindent{\em Case 2.2.} If there is $j_0\in\{1,\ldots, p\}$ such that $r_{j_0}>0$, then
\vsp{-0.2}
\begin{center}
$(g-1)(\textbf{x})=r_my_m+\cdots +r_{p+1}y_{p+1}+r_{p}+\cdots+(r_{j_0}-1)+\cdots +r_1$
\end{center}
\vsp{-0.2}
and $r_{j_0}-1\in [-1,0)$, so the induction hypothesis applies to $g-1$. Then there exists a formula $\varphi_3$ such that
$\varrho\circ (g-1)=\widetilde{\varphi_3}$.\\
\noindent{\em Case 2.3.} If there is $j_0\in\{p+1,\cdots, m\}$ such that $r_{j_0}>0$, then we set  $h_0(\textbf{x})=r_{j_0}y_{j_0}$ and
\vsp{-0.2}
\begin{center}
$g_0(\textbf{x})=g(\textbf{x})-r_{j_0}y_{j_0}-1.$
\end{center}
\vsp{-0.2}
It  follows that $g-1=g_0+h_0$ such that $g_0$ satisfies the induction hypothesis and $h_0:[0,1]^n\to [0,1]$. We are in the hypothesis of {\em Case 1}, so
there exists a formula $\varphi_3$ such that $\varrho\circ(g-1)=\widetilde{\varphi_3}$.\\
Summing up, we get  $\varrho\circ (g+h)  =\widetilde{\varphi}$ with $t={((\varphi_2\oplus \varphi_3^*)\odot \varphi_1)}$.
\end{proof}

\section{ Connections with the Pierce-Birkhoff conjecture}\label{pbc}

\noindent At the end of the paper \cite{BP}, the authors asked for a characterization of  the "{\em free, commutative, real $\ell$-algebra ($\ell$-group) with $n$ generators}" and they  conjectured that "{\em  it is isomorphic with the l-group of real functions which are continuous and piecewise polynomial of degree at most $n$ over a finite number of pieces}".  They asked "{\em the same problem
 for the free (commutative) $\ell$-rings, for free $f$-rings}", saying that: {\em "The former is probably very difficult"}.

\bdfn
Let $n\geq1$ be a natural number.
\begin{itemize}
\item A function  $f:{\mathbb R}^n\to {\mathbb R}$  is a  {\em piecewise polynomial} (PWP) function  if  it is continuous and there is a finite set of  polynomials   $\{p_1, \ldots, p_k\}\in {\mathbb R}[x_1,\ldots, x_n]$ such that
  for any $(a_1,\ldots, a_n)\in {\mathbb R}^n$  there exists  $i\in \{1,\ldots, k\}$ with
$f(a_1,\ldots, a_n)=p_i(a_1,\ldots, a_n)$.

\item   A continuous function  $f:{\mathbb R}^n\to {\mathbb R}$  is a  {\em inf-sup-polynomial-definable} (ISD) function  if  there is a finite set of  polynomials   $\{q_{ij}| 1\leq i\le m, 1\leq j\leq k\}\subseteq {\mathbb R}[x_1,\ldots, x_n]$ such that 
$f=\bigvee_{i=1}^m\bigwedge_{j=1}^k q_{ij}$.
\end{itemize}
\edfn

We denote by $PWP(n)$ the set of all PWP-functions and by $ISD(n)$ the set of all ISD-functions defined as above.

\bfact {\em The Pierce-Birkhoff conjecture} states that $ PWP(n)=ISD(n)$ for any $n\geq 2$ and, in this form, it was formulated by Henriksen and Isbell. The proof for $n\le2$ was made by L. Mah\'{e} in \cite{Mahe}, where  an unpublished proof of Gus Efroymson is also quoted. 
\efact

\bdfn
Let $n\geq 1$ be a natural number.

\begin{itemize}
\item A function  {$f: [0,1]^n\to [0,1]$}  is a   PWP$_u${\em-function}  if  it is continuous and there is a finite set of  polynomials   with real coefficients $p_1, \ldots, p_k: {\mathbb R}^n\to {\mathbb R}$ such that
  for any $(a_1,\ldots, a_n)\in {\mathbb R}^n$  there exists  $i\in \{1,\ldots, k\}$ with
$f(a_1,\ldots, a_n)=p_i(a_1,\ldots, a_n)$.

\item   A continuous function {$f: [0,1]^n\to [0,1]$}  is an   ISD$_u${\em-function}  if  there is a finite set of  polynomials with real coefficients  $\{q_{ij} :[0,1]^n\to {\mathbb R}| 1\leq i\le m, 1\leq j\leq k\}$ such that 
{$f=\bigvee_{i=1}^m\bigwedge_{j=1}^k ((q_{ij}\vee {\bf 0})\wedge {\bf 1})$}.
\end{itemize}
\edfn

We denote by $PWP(n)_u$ the set of all PWP$_u$-functions and by $ISD(n)_u$ the set of all ISD$_u$-functions.

\bthm \label{thmain}
The following properties hold:\\
 (1) $ISD(n)_u \subseteq FR_n\subseteq PWP(n)_u$ for any $n\in \mathbb{N}$,\\
 (2) $ISD(n)_u = FR_n =  PWP(n)_u$, for $n\le 2$.
\ethm
\begin{proof}
$(1)$ It is a direct consequence of Propositions \ref{formfct} and  \ref{ppf}.\\
$(2)$ Let $f\in PWP(n)_u$ with $n\le 2$. $f$ can be extended to a function $\overline{f}: \mathbb{R}^n \rightarrow \mathbb{R}$ in the following way: for $n=1$, we set $\overline{f}(x)=c_0$ for any $x\le 0$ and $\overline{f}(x)=c_1$ for any $x\ge 1$, where $c_0=f(0)$ and $c_1=f(1)$; for $n=2$, since $f$ is defined over $[0,1]^2$, the result follows by \cite[Theorem 1.2]{FM}. Therefore, by \cite{Mahe} there exist two set of indexes $I,\ J$ such that
$$\overline{f}= \bigvee _{i\in I} \bigwedge _{j\in J} \overline{f_{ij}},$$
where $\overline{f_{ij}}$ are polynomial functions from $\mathbb{R}^n$ to $\mathbb{R}$ with real coefficient. We consider the restrictions $f$ and $f_{ij}$ on $[0,1]^n$ of $\overline{f}$ and $\overline{f_{ij}}$ respectively, and we get \\
$$f=\varrho \circ f= \bigvee _{i\in I} \bigwedge _{j\in J} (\varrho \circ f_{ij})= \bigvee _{i\in I} \bigwedge _{j\in J} \widetilde{\varphi_{ij}},$$
for some suitable formulas, by Proposition \ref{formfct} . We set $\varphi=\bigvee _{i\in I} \bigwedge _{j\in J} \varphi_{ij}$, and we get $f=\widetilde{\varphi}$. Then $PWP(n)_u=FR_n=ISD(n)_u$, for $n\le 2$.
\end{proof}

\bfact We worked in the context of \textit{f}MV-algebras, so  the components of the piecewise polynomial functions have real coefficients. A similar approach can be used in the context of PMV-algebras, but in this case the components of the piecewise
polynomial functions will have integer coefficients. One can easily see that Propositions \ref{formfct}, \ref{ppf} and Theorem \ref{thmain}
can be proved in  a similar way. 
\efact

\begin{conj}\label{end}  $ISD(n)_u=FR_n=PWP(n)_u$ for any $n\geq 3$.
\end{conj}

\noindent{\bf Conclusion.}  The Pierce-Birkhoff conjecture implies Conjecture \ref{end} in the presence of additional extension results
(see \cite{FM} for $n=2$). Assuming Conjecture \ref{end}  holds, in order to prove the Pierce-Birkhoff conjecture, one needs results of the following type:

 {\em 
if $f\in PWP(n)$ such that $(f\wedge 1 \vee 0)|_{[0,1]^n}\in ISD(n)_u$ then $f\in ISD(n)$}.\\
 Conjecture \ref{end} is a   normal form theorem for the logical system  $\mathcal{FMVL}^{+}$, obtained by extending \L ukasiewicz logic with both an internal and external product and it can be seen  as a {\em local version} of the Pierce-Birkhoff conjecture.

\section*{Acknowledgements}
\noindent   Part of this research was carried on while S. Lapenta  was visiting University of Bucharest, supported by the grant of the International Doctoral Seminar School Entitled  "Janos Bolyai". Both authors  wish to thank Professor Antonio Di Nola for encouraging  their collaboration,  for his comments and discussions on  various subjects that, all together, led to the present development.

\bibliographystyle{elsarticle-num}

\end{document}